\documentclass[a4paper]{article}

\begin{document}

\title{ORTHOGONAL POLYNOMIALS IN ANALYTICAL METHOD OF SOLVING
DIFFERENTIAL EQUATIONS DESCRIBING DYNAMICS OF MULTILEVEL SYSTEMS\footnote{
{\it Integral Transforms and Special Functions}, 2000, Vol.\,10, No.\,3--4,
pp.\,299--308. Received December 12, 1999.\newline \copyright\, OPA (Overseas
Publishers Association) N.V. with permission from Gordon and Breach
Publishers.}}
\author{
V.A.~SAVVA${}^1$, V.I.~ZELENKOV${}^2$ and A.S.~MAZURENKO${}^1$\protect\\
${}^1$ {\it Institute of Physics, National Academy of Science, 68 Skorina
pr.,}\protect\\ {\it 220072 Minsk, Republic of Belarus}\protect\\
${}^2$ {\it International Sakharov Environmental University, 23
Dolgobrodskaya ul.,} \protect\\ {\it 220009 Minsk, Republic of Belarus} }
\date{}
\maketitle

\begin{abstract}
An effective method to obtain exact analytical solutions of equations
describing the coherent dynamics of multilevel systems is presented. The
method is based on the usage of orthogonal polynomials, integral transforms
and their discrete analogues. All the obtained solutions are expressed by
way of special or elementary functions.

Key words: integral transforms, orthogonal polynomials, multilevel quantum
systems.

1991 AMS Subject Classification: 33C45, 34C25
\end{abstract}

\section{INTRODUCTION}
Orthogonal polynomials and functions are used as the complete bases of the
solutions of differential and difference equations. In quantum physics
they are of primary importance in the construction of wave functions which
are orthogonal by definition. In present paper orthogonal polynomials and
functions are used to describe excitation of quantum multilevel systems in
high-power laser field.

The Schr\"odinger equation for a multilevel system being exciting in a
monochromatic laser field looks like
\begin{equation}\label{L1}
-id\/a_{n} (t)/d\/t=f_{n+1} \exp \left(-i\varepsilon_{n+1} t\right)
\,a_{n+1} (t)+f_{n} \exp \left(i\varepsilon_{n} t\right) \,a_{n-1} (t),\quad
n=0,\/1,\/2,\ldots
\end{equation}
where $f_{n}$ is a dipole moment function for the transitions
$n{-}1{\leftrightarrow}n$ (that is $\mu_{n-1,n}=f_{n} \mu_{0,1} $),
$\varepsilon_{n} $ is a frequency detuning. Probability amplitudes
$a_{n}(t)$ allow to determine levels populations
$\rho_n(t)=\left|a_n(t)\right|^2$, an average number of absorbed photons
$<n>=\sum_{n=1}^{n_{max}} n\rho_n(t)$ etc.

To obtain the solution of  (\ref{L1}) we introduce an analytical method
\cite{JMS} based on orthogonal functions $p_{n}(z)$.

\section{FROM ORTHOGONAL POLYNOMIALS TO QUANTUM SYSTEMS}
We seek the solution of  (\ref{L1}) in form
\begin{equation}\label{L2}
a_{n} (t)=\int_A^B \sigma(x) \frac{p_{0} (x)}{d_{0}} \frac{p_{n} (x)}{d_{n}}
\exp \left[it\left(r\theta (x)+s_{n} \right) \right] dx,
\end{equation}
where $d_{n}$ is the squared norm and $s_{n}$ determines the frequency
detuning: $\varepsilon_{n}=s_{n}-s_{n-1}$.

At $t=0$ the solution (\ref{L2}) passes into the orthogonality relation for
the functions $p_{n} (x)$. Hence we obtain the initial condition for Eq.\
(\ref{L1}): $a_{n}(t)=\delta_{n,0}$. The functions $p_{n} (x)$ should
satisfy the recurrence relation
\begin{equation}\label{L3}
\left(r\theta (x)+s_{n} \right) \frac{p_{n} (x)}{d_{n}}=f_{n}
\frac{p_{n-1} (x)}{d_{n-1}} +f_{n+1} \frac{p_{n+1} (x)}{d_{n+1}}.
\end{equation}

Initially we select the sequence of orthogonal functions which determines the
quantum system to be investigated. The takeoff criterion is the relevant
recurrence relation which allows to find both the constant $r$ and the
functions $f_{n}$ and $s_{n}$. After that we obtain the analytical solution
of  (\ref{L1}) using the integral transform like (\ref{L2}).

The application of analytical method (\ref{L1}--\ref{L3}) makes possible
describing dynamics of manifold quantum models including: $N$-level and
infinite-level, equidistant and non-equidistant, excited resonantly and
nonresonantly, closed and nonclosed, excited by monochromatic field and
laser pulse. Gegenbauer, Jacobi, Pollaczek, Laguerre, Meixner, Krawtchouk,
Hahn polynomials; Szeg{\"o}{-}Jacobi polynomials built using the Szeg{\"o}
method; Christoffel-Legendre orthogonal polynomials built by using the
Christoffel formula; orthogonal Gauss hypergeometric functions are used by
us to obtain the exact solutions for these systems.

This quantum systems manifold includes the systems with various
characteristics, first of all with qualitatively different dipole moment
functions $f_n$. In particular Jacobi, Legendre and Chebyshev of the
$1^{\mbox{\small st}}$ kind polynomials lead to the quantum systems with
$f_{n}\le 1$, Jacobi, Gegenbauer and Chebyshev of the $2^{\mbox{\small
nd}}$ kind polynomials~--- to the systems with $1\leq f_{n}<\sqrt{n}$,
Pollaczek, Laguerre and Meixner polynomials give rise to two sets of the
systems with $\sqrt{n}<f_{n}\leq n$ and $f_{n}>n$. It is known that
Hermite polynomials correspond to the harmonic oscillator ($f_n=\sqrt{n}$).

Thus the existing diversity of orthogonal polynomials makes possible to obtain
the solutions of  (\ref{L1}) describing the dynamics of various quantum
systems with different properties. Each orthogonal polynomials sequence (OPS)
can be put in unique correspondence with some quantum system. Below we shall
discuss some examples.

\section{JACOBI POLYNOMIALS}
Let's consider the quantum systems which have been built by means of
Jacobi polynomials. The dipole moment function following from a recurrence
relation (\ref{L3}) looks like
\begin{eqnarray}\label{V1}
f_{n} &= & \frac{2r}{2n+\alpha +\beta} \left[\frac{n\left(n+\alpha \right)
\left(n+\beta \right) \left(n+\alpha +\beta \right)}{\left( 2n+\alpha
+\beta -1\right) \left(2n+\alpha +\beta +1\right)} \right]^{^1/_2}, \\
r &= & \frac{\alpha +\beta +2}{2} \left[\frac{\alpha +\beta +3}{\left(
\alpha +1\right) \left(\beta +1\right)} \right]^{^1/_2},\quad \alpha
>-1,\ \beta >1.
\end{eqnarray}
Thus the dipole moments of higher transitions tend to the constant value
$r/2$. Energy levels of the Jacobi systems are nonequidistant but the
frequency detuning decreases with $n$ approximately as $n^{-3}$, i.e. the
highest levels are practically equidistant and they are being excited
resonantly.

Jacobi OPS contains two parameters $\alpha$ and $\beta$ which determine the
energy levels positions. Therefore we obtain two{-}parameter continual
family of multilevel quantum systems. Each pair of $(\alpha,\beta)$
corresponds to the specific Jacobi quantum system which is described by
dipole moment $f_n$ and frequency detuning $\varepsilon_n$ functions.

The solution of  (\ref{L1}) for the Jacobi systems dynamics is expressed
through the Kummer hypergeometric functions:
\begin{eqnarray}\label{V2}
a_{n}(t)=\left(2irt\right)^{n} \left[\frac{\left(\alpha +1\right)_{n}
\left(\beta +1\right)_{n}}{\left(\alpha +\beta +2\right)_{2n} \left(
n+\alpha +\beta +1\right)_{n} n!} \right]^{^1/_2} \times \nonumber\\
 \\
e^{it\left(r+s_{n} \right)} \/_{1} F_{1} \left(n+\alpha +1;2n+\alpha
+\beta +2;-2irt\right).\nonumber
\end{eqnarray}

At $\alpha=\beta$ Jacobi polynomials $P^{\left(\alpha,\alpha \right)} (x)$
pass into Gegenbauer ones $C^{\left(\lambda \right)}(x)$, $\lambda =\alpha
+1/2$. Thus Kummer function in  (\ref{V2}) is reduced to a cylindrical
function $J_{n+\lambda} (rt)$. The case $\alpha=\beta=1$ (Chebyshev
polynomials of the $2^{\mbox{\small nd}}$ kind) gives so called equal-Rabi
system $f_{n}\equiv 1$, for $\alpha=\beta=0$ (Chebyshev polynomials of the
$1^{\mbox{\small st}}$ kind) we obtain the system with
$f_{1}=1,$$f_{n}=1/\sqrt{2}$, $n\geq 2$. The levels of both Gegenbauer and
Chebyshev systems are equidistant, the excitation is resonant.

An interesting situation arises in a case $\alpha=-\beta$ when
$$
f_{n}=\left[3\left(n^{2} -\alpha^{2}\right) \bigm/ \left(1-\alpha^{2}
\right) \left(4n^{2} -1\right) \right]^{^1/_2},\quad s_{n}=\alpha
\left[3/\left(1-\alpha^{2} \right) \right]^{^1/_2} \delta_{n,0}.
$$
Here the frequency detuning $\varepsilon_{n}=-s_{0}\delta_{n,1}$, the
levels are equidistant everywhere except for the first transition and only
this transition is nonresonant. In particular at $\alpha=-\beta=\pm 1/2 $
we obtain $f_{n} \equiv 1$, $\varepsilon_{n} =-\delta_{n,1}$ and the
expressions for levels populations $\rho_{n} (t)=\left| \/a_{n}
(t)\right|^{2}$ reshape an especially simple form $\rho_{n}
\left(t\right)=J_{n}^{2} \left(2t\right)+ J_{n+1}^{2}\left(2t\right)$.

\section{KRAWTCHOUK POLYNOMIALS}
Using Krawtchouk polynomials of discrete variable \cite{NSU} we can
describe the nonresonant excitation of the quantum systems with finite
number of equidistant levels and dipole moment function
\begin{equation}\label{K1}
 f_n=\left[\frac{n\left(N-n+1\right)}{N}\right]^{^1/_2},\qquad
 n=\overline{0,N}, \quad N=1,2,\ldots.
\end{equation}

The analytical solution of  (\ref{K1}) is expressed \cite{Kraw} by way of
elementary functions:
\begin{eqnarray}\label{K2}
 a_n(t)&=&\left[ {N\choose n}y^n(t)\right]^{^1/_2} \left\{
\left[1-(1+\sigma)y(t)\right]^{^1/_2}+
 i\left[\sigma y(t)\right]^{^1/_2}\right\}^{N-n} \exp\left[ it\varepsilon(n-N/2)
 \right];\\
 \rho_n(t)&=& {N\choose n}y^n(t) \left[ 1-y(t) \right]^{N-n};\\
 <n>&=&Ny(t)
\end{eqnarray}
where
\begin{equation}\label{K3}
 y\left(t\right)=\frac {\sin^2\left\{ \left[ (1+\sigma)\left/ N \right.
 \right]^{^1/_2}t \right\}} {(1+\sigma)} , \qquad
 \sigma=N\varepsilon^2/4.
\end{equation}

The Krawtchouk systems are suitable exactly solvable models for
vibrational excitation of a molecule in high-power infrared laser fields.
The solution describes an excitation with arbitrary detuning of laser
frequency from the system transitions frequency. The case $p=^1/_2$
corresponds to the resonant excitation. Krawtchouk systems family includes
both two-level and three-level systems and harmonic oscillator as special
and limit cases ($N=1$, $N=2$ and $N\to\infty$ correspondingly).

Thus the Krawtchouk polynomials give a possibility to build one-parametric
family of quantum multilevel systems and to obtain exact analytical
solutions of the equations describing the excitation dynamics. The
Krawtchouk systems family is an extensive class of multilevel models for
various dynamical processes in spectroscopy, nonlinear optics and the
other fields.

The method under consideration gives also the possibility of solving the
equations which are more complicated then  (\ref{L1}). Let the quantum system
consist of finite number of degenerate levels. The equation describing the
resonant excitation of such a system is
\begin{eqnarray}
&&\!\!\!-i\frac{da_{n,m}(t)}{dt} =\Omega
\left[f_{n+1}a_{n+1,m}(t) +
f_na_{n-1,m}(t)\right] +
\label{K4} \\
&&\!\!\!\Omega^{\prime}\left[f_{n+1}g_{m+1}a_{n+1,m+1}(t) +
f_{n+1}g_ma_{n+1,m-1}(t)\right. +
\vphantom{\frac{da_{n,m}(t)}{dt}}\left.
f_ng_{m+1}a_{n-1,m+1}(t)+f_ng_ma_{n-1,m-1}(t)\right] \nonumber
\end{eqnarray}
where $\Omega$ and $\Omega^{\prime}$ are Rabi frequencies, $f_n=\left[
n(N-n+1)/N\right]^{^1/_2}$ and $g_m=\left[ m(M-m+1)/M\right] ^{^1/_2}$ are
dipole moment functions for the interlevel and intersublevel transitions
correspondingly, $n=\overline{0,N}$, $m=\overline{0,M}$; $M,N=1,2,\ldots$
Here we assume $a_{n,m}(0)=\delta_{n,0} \delta_{m,0}$ as an initial
condition and $\Delta n=\pm 1$, $\Delta m=0,\pm 1$ as the selection rules.

The solution of  (\ref{K4}) is expressed via Krawtchouk polynomials and
elementary functions:
\begin{equation}
a_{n,m}=\frac{i^n}{2^{M-m}}\left[\frac{{N\choose n}}{{M\choose
m}}\right]^{^1/_2}
\sum\limits_{j=0}^Mk_m^{(^1/_2)}(j,M)\,{M\choose j}\sin^n\tau_j\,\cos
^{N-n}\tau_j \label{K5}
\end{equation}
where $\tau_j=\left({t}/{}_{\sqrt{N}}\right) \left[ \Omega +\Omega^{\prime}
(2j-M)\left/\sqrt{M}\right.\right]$.

Selecting other orthogonal polynomials we can obtain the solutions of Eq.\
(\ref{K4}) for different types of $f_n$ and $g_m$.

\section{CHRISTOFFEL-LEGENDRE POLYNOMIALS}
Let's consider the solution of  (\ref{L1}) using Christoffel-Legendre
orthogonal polynomials \cite{JAS}
\begin{equation}\label{L4}
p_{n}(x)=\frac{b}{b-x} \left[P_{n} (x)-\kappa_{n}^{-1}
(b)P_{n+1}(x)\right],\qquad \kappa_{n}(b)=P_{n+1}(b)/P_{n}(b)
\end{equation}
built with Christoffel formula \cite{Szego} using Legendre polynomials
$P_{n}(x)$. The polynomials $p_{n} (x)$ are orthogonal on an interval
$\left[-1,1\right]$ with respect to the weight
\begin{equation}\label{L5}
 \sigma (x)=(b-x)/b,\quad \quad \quad b\geq 1
\end{equation}
and squared norm
\begin{equation}\label{L6}
 d_{n}=\left\{\,\left[2/\left(n+1\right) \right] b/\kappa_{n} \left(
b\right) \right\}^{^1/_2}
\end{equation}
They satisfy the recurrence relation (\ref{L3}) ($\theta (x)\equiv x$)
where
\begin{eqnarray}\label{L7}
f_{n}&=&r\left(n/(2n+1) \right) d_{n-1} /d_{n},\nonumber\\
r&=&3d_{1} /d_{0},\\
s_{n}&=&r\left(\frac{n+1}{2n+1} \kappa_{n}
(b)-\frac{n+2}{2n+3} \kappa_{n+1} (b)\right).\nonumber
\end{eqnarray}

The application of Christoffel-Legendre polynomials in  (\ref{L2}) gives
\cite{JAS} the analytical solution
\begin{equation}\label{L8}
 a_{n} (t,b)=i^{n} \left(\frac{\pi}{2} \right)^{\frac{1}{2}}
\frac{d_{0}}{d_{n}} \,\frac{\exp \left(is_{n} \,t\right)}{(r\,t)^{^1/_2}}
\left[J_{n+0,5} \left(r\,t\right) -i\,\kappa_{n}^{-1} \left(b\right)
\,J_{n+1,5} \left(r\,t\right) \right]
\end{equation}
of  (\ref{L1}). This solution describes the excitation of one-parameter family
of multilevel Christoffel-Legendre systems with nonequidistant energy levels,
i.e. $\varepsilon_{n} \neq \mbox{const}$.

\section{LEGENDRE FUNCTIONS OF THE FIRST KIND}
Another quantum system can be built by using the functions
\begin{equation}\label{L9}
 p_{n}^{\lambda} (x)\,=\,\left. (-i)^{\mu} \Gamma (\mu
+1)P_{-\lambda}^{-\lambda -\mu} (x)\right|_{\mu=n},
\end{equation}
being a special case of Legendre function of the $1^{\mbox{\small st}}$
kind $P_{\nu}^{\mu} (x)$. The functions $p_{n}^{\lambda} (x)$ satisfy the
recurrence relation (\ref{L3}) where
\begin{eqnarray}\label{L10}
f_{n}&=&\left| \frac{n\left(2\lambda +n-1\right)\,(\lambda
+1)}{2(\lambda +n-1)\,(\lambda +n)} \right|^{^1/_2},\qquad
s_{n}\equiv 0, \nonumber\\
 \\
r&=&i\sqrt{2\left| \lambda +1\right|},\qquad
\theta(x)=\frac{x}{\sqrt{x^{2} -1}},\nonumber
\end{eqnarray}
and
\begin{equation}\label{L11}
d_{n}=\left| \vphantom{p^{I}_{I}}\Gamma (1-2\lambda -n)\Gamma (n+1)\left/
\left(\lambda +n\right)\right.\right|^{^1/_2},\quad \lambda <-1, \quad
0\leq n\leq \left|\lambda\right|.
\end{equation}

For $p_{n}^{\lambda}(x)$ is correct
\begin{equation}\label{L12}
p_{n}^{\lambda} (x)\,=\,2^{\lambda} i^{n} \,\frac{\Gamma \,(-2\lambda
-n+1)n!}{\Gamma \,\left(-\lambda +1\right)} \left(x^{2}
-1\right)^{-\lambda /2} C_{n}^{\lambda} \left(\frac{x}{\sqrt{x^{2} -1}}
\right),\quad \lambda <-1,\quad 0\leq n\leq \,\,|\lambda |
\end{equation}
where $C_{n}^{\lambda}$ are Gegenbauer functions \cite{BE} of a complex
argument $z$
\begin{equation}\label{L13}
C_{\alpha}^{\lambda} (z)=\frac{\Gamma (\alpha +2\lambda)}{\Gamma (\alpha
+1)\Gamma (2\lambda)} {}_2\mbox{\Large $F$}_1 \left( {{-\alpha,\;
\alpha+2\lambda} \atop {\lambda+1/2}} \left| \frac{1-z}{2} \right. \right).
\end{equation}

Functions $p_{n}^{\lambda}$ are orthogonal on the discrete interval:
\begin{equation}\label{L14}
\sum\limits_{l=0}^{N-1}\sigma^{(N)}(x_l) \/p_{m}^{\lambda}(x_l)
p_{n}^{\lambda}(x_l)=\delta_{m,\,n} d_{n}^{2},\quad N>\left| 1-\lambda
\right|,\quad \lambda<-1
\end{equation}
with respect to the weight function
\begin{equation}\label{L15}
\sigma^{(N)}(x_l)= \frac{{\cal K}_{N}}{{\cal K}_{N-1}}
\left(\frac{2^{-\lambda} \Gamma (1-\lambda)d_{N-1}}{\Gamma (N)\Gamma
(2-2\lambda -N)} \right)^{2} \left(x_l^{2} -1\right)^{\lambda} \left(\left.
C_{N-1}^{\lambda }(z_l) \frac{d}{d\,z_l} C_{N}^{\lambda}(z_l)\right|_{z_l=
x_l\left/ \sqrt{x_{l}^{2}-1}\right.} \right)^{-1}
\end{equation}
where $z_l$ are the roots of a polynomial of degree $N$, $
C_N^{\lambda}(z_l)=0, x_{l}=z_{l}/ \sqrt{z_{l}^{2}-1}$ and ${\cal K}_{N}$
is the higher coefficient of this polynomial. Thus the solution of
(\ref{L1}) for the system of $N$ equidistant levels in a resonant field
(at $f_{N} \equiv 0$ and $\varepsilon_{n} \equiv 0$) looks like
\begin{equation}\label{L16}
a_{n}^{(N)} (t)=\sum\limits_{l=0}^{N-1}\,\sigma^{(N)}(x_{l})\/\left(
\frac{p_{m}^{\lambda} (x_{l})}{d_{m}} \right)^{\!\!*} \frac{p_{n}^{\lambda}
(x_{l})}{d_{n}} \exp \,\left[\frac{itr\,x_{l}}{\sqrt{x_{l}^{2} -1} }
\right].
\end{equation}

\section{FROM QUANTUM SYSTEMS TO ORTHOGONAL POLYNOMIALS}
Previously we initially selected orthogonal functions to obtain the solution
of the equation (\ref{L1}) for qualitatively and quantitatively different
quantum systems. Having done so we get the compact solutions which are
expressed via known special or even elementary functions. When OPS contains
the parameters the solution describes the dynamics of continual or discrete
family of quantum systems. From the other side the imperfection of this method
is that it is not always possible to obtain exact analytical solutions for a
target system with the given dipole moments function $f_n$.

Let's consider the quantum system with an arbitrary dipole moments function
being restricted with the case of resonant excitation and a finite number of
equidistant levels. Now  (\ref{L1}) reduces to the system of $N+1$ uniform
linear differential equations
\begin{eqnarray}\label{V3}
-i\ d\/a_{0} (t)/d\/t &= & f_{1} a_{1} (t); \nonumber\\
-i\ d\/a_{n} (t)/d\/t &= & f_{n+1} a_{n+1} (t)+f_{n} a_{n-1} (t);
\quad n=\overline{1,N-1} \\
-i\ d\/a_{N} (t)/d\/t &= & f_{N} a_{N-1} (t)\nonumber
\end{eqnarray}
(some special cases of $f_{n}$ have been studied in \cite{BB}).

Let us seek the solution of  (\ref{V3}) in the form
\begin{equation}\label{V4}
a_{n} (t)=\sum\limits_{k=0}^{N}\sigma_{k} p_{n} \left(\Lambda_{k}
\right) \exp \left(i\Lambda_{k} t\right)
\end{equation}
where $p_{n}$, $\Lambda_{k}$ and $\sigma_{k}$ are the functions to be
determined. Substituting (\ref{V4}) into (\ref{V3}) we come to a three-term
recurrence relation for $p_{n}$:
\begin{equation}\label{V5}
f_{n+1} p_{n+1} \left(\Lambda_{k} \right) +f_{n} p_{n-1} \left(\Lambda
_{k} \right)=\Lambda_{k} p_{n} \left(\Lambda_{k} \right),\quad p_{-1}=0,
\quad p_{0}=\mbox{const}
\end{equation}
Therefore $p_{n}$ are the polynomials orthogonal with respect to variable
$\Lambda_{k}$ (generally on a nonuniform grid).

This truncated orthogonal polynomials $p_{n} \left(\Lambda_{k} \right)$
differ from the common polynomials listed above which satisfy the same
dipole moment function $f_{n}$ with an interval of $n$ change. Let's
discover the basic characteristics of these polynomials.

It is known that the recurrence relation like (\ref{L3}) determines an
orthogonal polynomial up to an arbitrary constant factor. The comparison of
 (\ref{L3}) and  (\ref{V5}) shows that the polynomials $p_{n}
\left(\Lambda_{k} \right)$ are orthonormalized. Therefore it is natural to
choose the normalization $p_{0}\left(\Lambda_{k} \right)\equiv 1$.

Suppose that initially all the particles are concentrated on the lower level,
in other words $a_{n} \left(0\right)=\delta_{n,0}$. Then it follows from
(\ref{V4}) that
\begin{equation}\label{V6}
\sum\limits_{k=0}^{N}\sigma_{k} p_{n} \left(\Lambda_{k} \right) p_{0}
\left(\Lambda_{k} \right)=\delta_{n,0}.
\end{equation}
Thus the integration constants $\sigma_{k}$ are the weight coefficients of
polynomials $p_{n} $ with respect to the variable $\Lambda_{k}$.

To determine $\Lambda_{k} $ we shall use the last equation of (\ref{V3});
it gives
\begin{equation}\label{V7}
f_{N} p_{N-1} \left(\Lambda_{k} \right)=\Lambda_{k} p_{N}
\left(\Lambda_{k} \right).
\end{equation}
Comparing  (\ref{V7}) and  (\ref{V5}) it is possible to conclude that
$\Lambda_{k} $ are the roots of the polynomial $p_{N+1} \left(\Lambda_{k}
\right)$.

Calculation of the polynomials $p_{n} \left(\Lambda_{k} \right)$ by the
recurrence relation (\ref{V5}) and determination of weight coefficients
$\sigma_k$ from the system (\ref{V6}) are possible for any values of $N$. At
the same time  (\ref{V6}) is reduced eventually to the solving in radicals of
the polynomial equations of the degree $N{+}1$. Generally it is possible only
for $N{+}1\leq 4$. However in the case of a resonant excitation the
polynomials $p_{n} $ contain only even or odd degrees of an argument
(depending on $n$ parity). Thus the roots of the polynomials up to the ninth
degree can be always found analytically. However using some computer algebra
system (e.g.\,{\tt Mathematica}) it is possible to obtain the roots for the
polynomials of higher degrees for some simple $f_{n}$. In any case it is
possible to find these roots numerically with an arbitrary precision.

The method allows to obtain the solutions describing the dynamics of the
quantum systems with an arbitrary dipole moment dependence on level number.
When the solution of  (\ref{L1}) with the same $f_{n}$ but for infinite number
of levels can be received under the formula (\ref{L2}) by selecting the
relevant orthogonal polynomials it is interesting to compare the truncated
polynomials with the common ones. From the comparison of  (\ref{L3}) and Eq.\
(\ref{V5}) it follows that
\begin{equation}\label{V8}
p_{n} \left(\Lambda_{k} \right)=p_{n}^{(c)} \left(\Lambda_{k} /r \right)
\end{equation}
where $p_{n}^{(c)}$ are common polynomials, $r=b_{0}d_{0}/d_{1}$
\cite{Szego}, $d_{n}$ is a squared norm, $b_{0}$ is defined from the
recurrence relations for unnormalized polynomials: $p_{n+1} (x)+b_{n}
p_{n-1} (x)=$ $\left(rx+s_{n} \right) p_{n} (x)$. In particular if
$R_{k}^{\left( N+1\right)} $ is the $k^{\mbox{\small th}}$ root of
$p_{N+1}^{(c)} $ then $\Lambda_{k}=rR_{k}^{\left(N+1\right)}$.

Methods (\ref{L1}--\ref{L3}) and (\ref{V4}) can be generalized to the
system excitation in an arbitrary laser field, to the account of
non-adjacent dipole transitions and for the systems of degenerate levels.

\end{document}